\documentclass[11pt]{article}
\usepackage{psfig}
\usepackage{amsmath}
\usepackage[numbers, square, comma]{natbib} 
\usepackage{geometry}                
\usepackage{graphicx}
\usepackage{amssymb}
\usepackage{epstopdf}
\usepackage {natbib}

\begin{document}

\fontsize{11}{14.5pt}\selectfont

\begin{center}

{\small Technical Report No.\ 0510,
 Department of Statistics, University of Toronto}

\vspace*{0.9in}

{\Large\bf Improving Classification When a Class Hierarchy is Available \\[4pt] Using a Hierarchy-Based Prior}\\[16pt]

 \parbox[t]{7.6cm} {
 \begin{center}
{\large Babak Shahbaba}\\[2pt]
 Dept.\ of Public Health Sciences, Biostatistics\\ 
 University of Toronto \\
 Toronto, Ontario, Canada \\
 \texttt{babak@stat.utoronto.ca}\\[10pt]
 \end{center}
}
\hfill
\parbox[t]{8.5cm} {
\begin{center}
{\large Radford M. Neal}\\[2pt]
 Dept.\ of Statistics and Dept.\ of Computer Science\\
 University of Toronto \\
 Toronto, Ontario, Canada \\
 \texttt{radford@stat.utoronto.ca}\\[10pt]
 \end{center}
 }

18 October 2005

\end{center}

\vspace*{8pt}

\noindent \textbf{Abstract.}  We introduce a new method for building
classification models when we have prior knowledge of how the classes
can be arranged in a hierarchy, based on how easily they can be
distinguished. The new method uses a Bayesian form of the multinomial
logit (MNL, a.k.a. ``softmax'') model, with a prior that introduces
correlations between the parameters for classes that are nearby in the
tree. We compare the performance on simulated data of the new method,
the ordinary MNL model, and a model that uses the hierarchy in
different way.  We also test the new method on a document labelling
problem, and find that it performs better than the other methods,
particularly when the amount of training data is small.

\section{\hspace*{-7pt}Introduction}\label{sec-intro}\vspace*{-12pt}

In this paper, we consider classification problems where classes have a hierarchical structure. The hierarchy reflects our prior opinion regarding similarity of classes. Two classes are considered similar if it is difficult to distinguish them from each other on the basis of the features available. The similarity of classes increases as we descend the hierarchy. 

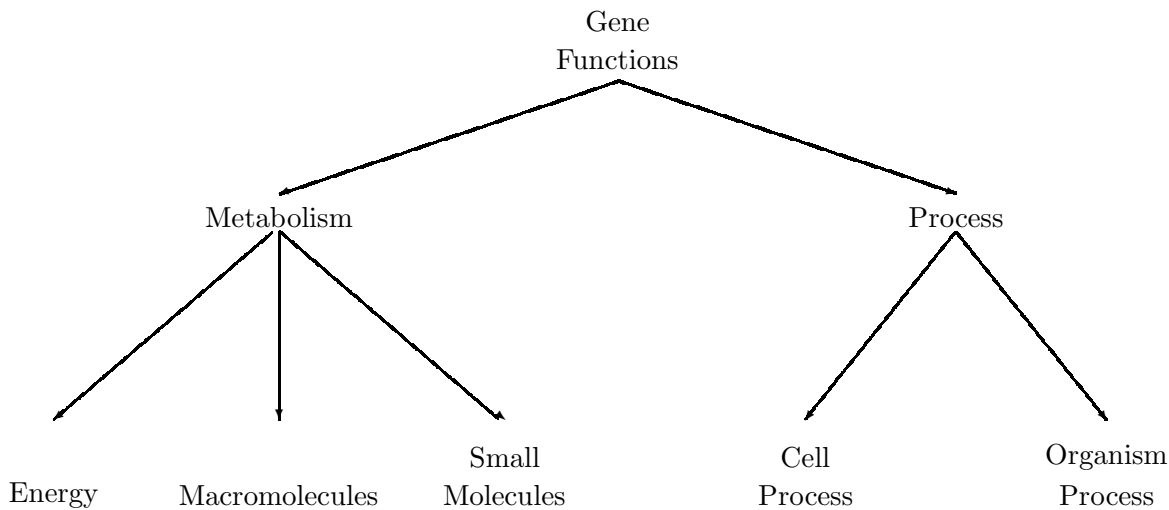
\begin{figure}[h] 
\begin{center}

\unitlength 1mm
\begin{picture}(150,70)(5,5)
\put(60,55){\makebox(0,0)[cc]{}}

\put(85,73){\makebox(0,0)[cc]{Gene}}

\put(85,68){\makebox(0,0)[cc]{Functions}}

\put(60,55){\makebox(0,0)[cc]{}}

\linethickness{0.3mm}
\multiput(40,50)(0.36,0.12){125}{\line(1,0){0.36}}
\put(40,50){\vector(-3,-1){0.12}}
\linethickness{0.3mm}
\multiput(85,65)(0.36,-0.12){125}{\line(1,0){0.36}}
\put(130,50){\vector(3,-1){0.12}}
\linethickness{0.3mm}
\multiput(10,20)(0.14,0.12){208}{\line(1,0){0.14}}
\put(10,20){\vector(-4,-3){0.12}}
\linethickness{0.3mm}
\put(40,20){\line(0,1){25}}
\put(40,20){\vector(0,-1){0.12}}
\linethickness{0.3mm}
\multiput(40,45)(0.14,-0.12){208}{\line(1,0){0.14}}
\put(70,20){\vector(4,-3){0.12}}
\linethickness{0.3mm}
\multiput(110,20)(0.12,0.15){167}{\line(0,1){0.15}}
\put(110,20){\vector(-3,-4){0.12}}
\linethickness{0.3mm}
\multiput(130,45)(0.12,-0.15){167}{\line(0,-1){0.15}}
\put(150,20){\vector(3,-4){0.12}}
\put(40,47){\makebox(0,0)[cc]{Metabolism}}

\put(130,47){\makebox(0,0)[cc]{Process}}

\put(10,10){\makebox(0,0)[cc]{Energy}}

\put(70,15){\makebox(0,0)[cc]{Small}}

\put(70,10){\makebox(0,0)[cc]{Molecules}}

\put(110,15){\makebox(0,0)[cc]{Cell}}

\put(110,10){\makebox(0,0)[cc]{Process}}

\put(150,15){\makebox(0,0)[cc]{Organism}}

\put(150,10){\makebox(0,0)[cc]{Process}}

\put(40,10){\makebox(0,0)[cc]{Macromolecules}}

\end{picture}

\caption{A part of a gene annotation hierarchy.} 
 \label{geneAnnotation} 
\end{center}
\end{figure}

One example of a hierarchical classification scheme is the annotations describing the biological function of genes. Functions are usually presented in a hierarchical form starting with very general classes and becoming more specific in lower levels of the hierarchy. There are many different gene annotation schemes. Figure \ref{geneAnnotation} shows a small part of the scheme proposed by Rison \cite{rison00}.

Another example of hierarchical classification is the document labelling problem. This involves classification of document regions to one of many predefined classes. It is possible to arrange these classes in a hierarchy such as that of Figure \ref{documentLabelling} below. We arranged these classes in a hierarchical form based on our guess as to how easily they can be distinguished. For instance, we considered ``Figure Caption'' and ``Table Caption'' similar to each other since they both usually contain small amounts of straight text and both are centre justified. 

Bayesian models provide a framework that allows us to incorporate prior knowledge of this sort. The prior distribution in Bayesian models represents our belief (as well as our uncertainty) regarding likely values of parameters. For modelling hierarchical classes, we introduce a new method which is based on a Bayesian form of the multinomial logit (MNL) model. We use a prior that introduces correlations between the parameters of nearby classes. 

This paper is organized as folllows. In section 2, simple classification models and their extensions for analysing hierarchical classes are discussed. In section 3, using simulated data, we compare the performance of our model, the ordinary MNL model and an alternative model that uses the hierarchy in a different way. In section 4 we compare the same models on a document labelling problem. The last section summarizes our findings and presents some ideas for future research.

\section{\hspace*{-7pt}Hierarchical classification}\label{sec-one}\vspace{-12pt}

Consider a classification problem in which we have observed data for $n$ cases, ($x^{(1)},y^{(1)}$), ...,($x^{(n)},y^{(n)}$), where $x^{(i)} = x_{1}^{(i)}, ..., x_{p}^{(i)}$ is the vector of $p$ covariates (features) for case $i$, and $y^{(i)}$ is the associated class. Our goal is to develop a model that assigns each data item to its correct class based on the observed covariates. The resulting model will be used to classify future cases for which the class membership is unknown but the covariates are available. For binary classification problems, a simple logistic model can be used:
\begin{eqnarray}
P(y=1|x, \alpha, \boldsymbol{ \beta}) & = & \frac{\exp(\alpha + x \boldsymbol{ \beta})}{1+\exp(\alpha+x \boldsymbol{\beta})}
\end{eqnarray}
Here, $\alpha$ is the intercept, $\boldsymbol{\beta}$ is a $p \times 1$ vector of unknown parameters and $x\boldsymbol{\beta}$ is its inner product with the covariate vector.

When there are three or more classes, we can use a generalization known as the multinomial logit (MNL) model: 
\begin{eqnarray} \label{mnl}
P(y=j|x, \boldsymbol{\alpha}, \boldsymbol{\beta}) & = &\frac{\exp(\alpha_j + x\boldsymbol{ \beta}_j)}{\sum_{j=1}^{c} \exp(\alpha_j+x\boldsymbol{ \beta}_j)}
\end{eqnarray}
where $c$ is the number of classes. For each class, $j$, there is a vector of $p$ unknown parameters $\boldsymbol \beta_j$. The entire set of regression coefficients $\boldsymbol{\beta = \beta_1, ..., \beta_c}$ can be presented as a $p \times c$ matrix. This representation is redundant, since one of the $\beta_j$'s can be set to zero without changing the set of relationships expressible with the model, but removing this redundancy would make it difficult to specify a prior that treats all classes symmetrically. For this model we can use the following priors:
\begin{eqnarray*}
\alpha_{j} | \tau_{0}   & \sim & N(0, \tau_{0}^{2}) \\
\beta_{jl} | \tau \; & \sim & N({0}, \tau^{2}) \\
\tau_{0}^{-2} & \sim & Gamma(a_{0}, b_{0})  \\
\tau^{-2} & \sim & Gamma(a, b)  \\
\end{eqnarray*}  
where $j = 1, ..., c$ and $l = 1, ..., p$.

The MNL model treats classes as unrelated entities without any hierarchical structure. This is not always a realistic assumption. In many classification problems, like those discussed above, one can arrange classes in a hierarchical form analogous to the hierarchy of species arranged in genera, families, etc. If the classes have in fact the assumed structure, one would expect to obtain a higher performance by using this additional information. A special case is when the classes are ordered (e.g., education level). For these problems a more parsimonious model (e.g., cumulative logit model) with improved power can be used \cite{agresti02}. 

One approach for modelling hierarchical classes is to decompose the classification model into nested models (e.g., logistic or MNL). Nested MNL models are extensively discussed in econometrics (e.g., \cite{sattath77, mcfadden81}) in the context of estimating the probability of a person choosing a specific alternative (i.e., class) from a discrete set of options (e.g., different modes of transportation). These models, known as discrete choice models, aim at forecasting and explaining human decisions through optimizing an assumed utility (preference) function, which is different from our aim of maximizing classification accuracy. 

Goodman \cite{goodman01} showed that using hierarchical classes can significantly reduce the training time of maximum entropy-based language models and results in slightly lower perplexities. He illustrated his approach using a word labelling problem, and recommended that instead of predicting words directly, we can first predict the class the word belongs to, and then predict the word itself. This approach can be applied to the other learning techniques.

\begin{figure} \label{tree}
\unitlength 1mm
\begin{picture}(147,75)(-2,0)
\put(55,60){\makebox(0,0)[cc]{$\boldsymbol{\phi_{11}}$ }}

\put(110,60){\makebox(0,0)[cc]{$\boldsymbol{\phi_{12}}$ }}

\put(31.5,36){\makebox(0,0)[cc]{}}

\put(19,28){\makebox(0,0)[cc]{$\boldsymbol{\phi_{21}}$ }}

\put(42,28){\makebox(0,0)[cc]{$\boldsymbol{\phi_{22}}$ }}

\put(119,28){\makebox(0,0)[cc]{$\boldsymbol{\phi_{23}}$ }}

\put(142,28){\makebox(0,0)[cc]{$\boldsymbol{\phi_{24}}$ }}

\put(15,9){\makebox(0,0)[cc]{Class 1}}

\put(45,9){\makebox(0,0)[cc]{Class 2}}

\put(115,9){\makebox(0,0)[cc]{Class 3}}

\put(145,9){\makebox(0,0)[cc]{Class 4}}

\put(13,2){\makebox(0,0)[cc]{$\boldsymbol{\beta_1 = \phi_{11} + \phi_{21} }$}}

\put(45,2){\makebox(0,0)[cc]{$\boldsymbol{\beta_2 = \phi_{11} + \phi_{22} }$}}

\put(115,2){\makebox(0,0)[cc]{$\boldsymbol{\beta_3 = \phi_{12} + \phi_{23} }$}}

\put(147,2){\makebox(0,0)[cc]{$\boldsymbol{\beta_4 = \phi_{12} + \phi_{34} }$}}

\linethickness{0.3mm}
\multiput(15,15)(0.12,0.2){125}{\line(0,1){0.2}}
\linethickness{0.3mm}
\multiput(30,40)(0.12,-0.2){125}{\line(0,-1){0.2}}
\linethickness{0.3mm}
\linethickness{0.3mm}
\multiput(115,15)(0.12,0.2){125}{\line(0,1){0.2}}
\linethickness{0.3mm}
\multiput(130,40)(0.12,-0.2){125}{\line(0,-1){0.2}}
\linethickness{0.3mm}
\put(80,70){\circle{8.33}}

\put(80,70){\makebox(0,0)[cc]{1}}

\put(105,75){\makebox(0,0)[cc]{}}

\linethickness{0.3mm}
\multiput(30,50)(0.4,0.12){125}{\line(1,0){0.4}}
\linethickness{0.3mm}
\multiput(80,65)(0.4,-0.12){125}{\line(1,0){0.4}}
\linethickness{0.3mm}
\put(30,45){\circle{8.33}}

\linethickness{0.3mm}
\put(130,45){\circle{8.33}}

\put(30,45){\makebox(0,0)[cc]{2}}

\put(130,45){\makebox(0,0)[cc]{3}}

\end{picture}

\caption{A simple representation of our model. The coefficient parameters for each classes are presented as a sum of parameters at different level of hierarchy }
\label{simpleHierarchy}
\end{figure}
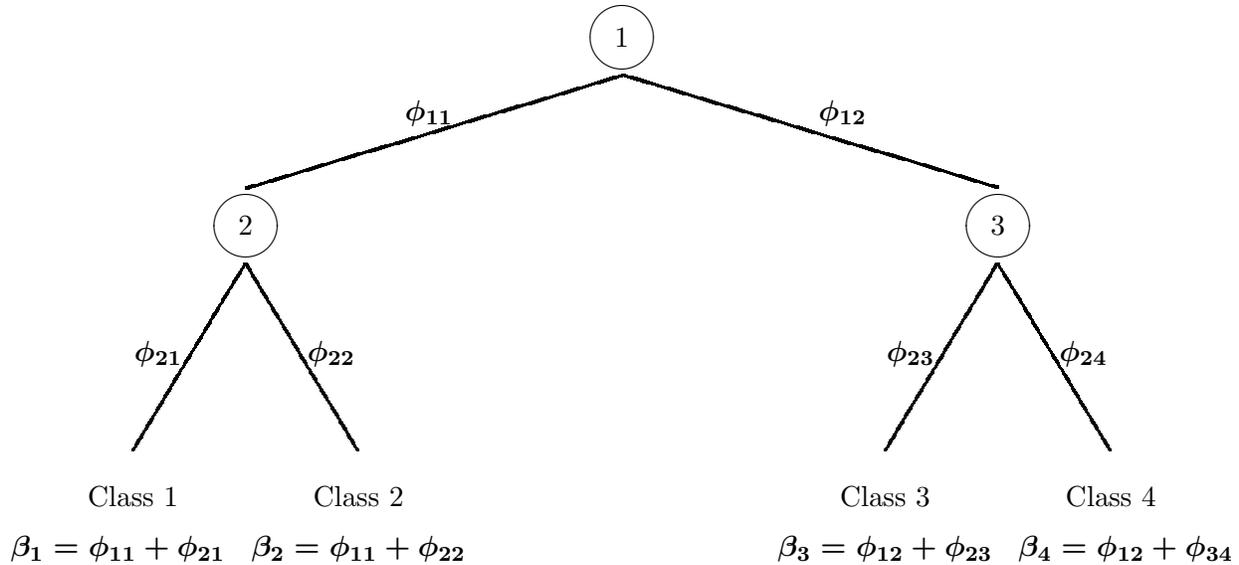

Fox \cite{fox97} suggested using successive binary partitions (i.e., dichotomies) of classes and fitting a logistic model to each partition. In other words, we recursively split a set of classes into two mutually exclusive subsets where each subset contains similar classes. In Figure \ref{simpleHierarchy}, these partitions are \{12, 34\}, \{1, 2\}, and \{3, 4\}. The resulting nested binary models are statistically independent. The likelihood can therefore be written as the product of the likelihoods for each of the binary models. For example, in Figure \ref{simpleHierarchy} we have
\begin{eqnarray}
P(y=1| x) & = & P(y \in \{1, 2 \} | x) \times P(y \in \{1\}| y \in \{1, 2\}, x)
\end{eqnarray}
Restriction to binary models is unnecessary. At each level, classes can be divided into more than two subsets and MNL can be used instead of logistic regression. We refer to methods based on decomposing the tree structure into nested MNL models as treeMNL. Consider a parent node, $m$, with $c_m$ children nodes. The portion of the nested MNL model for this node has the following form:
\begin{eqnarray*}
P(y=k|x, \boldsymbol{\alpha}_m, \boldsymbol{\beta}_m) & = &\frac{\exp(\alpha_{mk} + x\boldsymbol{ \beta}_{mk})}{\sum_{k=1}^{c_m} \exp(\alpha_{mk}+x\boldsymbol{ \beta}_{mk})} \\
\alpha_{mk} | \tau_{0} \  & \sim & N(0, \tau_{0}^{2}) \\
\beta_{mkl} | \tau_{m} & \sim & N({0}, \tau_{m}^{2}) \\
\tau_{0}^{-2} & \sim & Gamma(a_{0}, b_{0})  \\
\tau_{m}^{-2} & \sim & Gamma(a_m, b_m)  \\
\end{eqnarray*}  
where $k = 1, ..., c_m$ and $l = 1, ..., p$. We calculate the probability of each end node, $j$, by multiplying the probabilities of all intermediate nodes leading to $j$.

We introduce an alternative Bayesian framework for modelling hierarchical classes. Consider Figure \ref{simpleHierarchy}, which shows a hierarchical classification problem with four classes. For each branch in the hierarchy, we define a different set of parameters. In Figure \ref{simpleHierarchy}, these parameters are denoted as $\boldsymbol{\phi_{11}$ and $\phi_{12}}$ for branches in the first level and $\boldsymbol{\phi_{21}}$, $\boldsymbol{\phi_{22}}$, $\boldsymbol{\phi_{23}}$ and $\boldsymbol{\phi_{24}}$ for branches in the second level. We assign objects to one of the end nodes using a MNL model (equation \ref{mnl}) whose regression coefficients for class $j$ are represented by the sum of parameters on all the branches leading to that class. In Figure \ref{simpleHierarchy}, these coefficients are $\boldsymbol{\beta_{1} = \phi _{11}+\phi _{21} }$, $\boldsymbol{\beta_{2} = \phi_{11}+\phi _{22} }$, $\boldsymbol{\beta_{3} = \phi _{12}+\phi _{31} }$ and $\boldsymbol{\beta_{4} = \phi _{12}+\phi _{32} }$ for classes $1, 2, 3$ and $4$ respectively. Sharing the common terms, $\boldsymbol{\phi_{11}}$ and $\boldsymbol{\phi_{12}}$, introduces prior correlation between the parameters of nearby classes in the hierarchy.

In our model, henceforth called corMNL, $\boldsymbol{\phi}$'s are vectors with the same size as $\boldsymbol{\beta}$'s. We assume all the components of the $\boldsymbol\phi$'s are independent, and have normal prior distributions with zero mean. The variance of the $\phi$'s is regarded as a hyperparameter, which controls the magnitude of coefficient parameters. When a part of the hierarchy is irrelevant, we hope the posterior distribution of its corresponding hyperparameter will be concentrated near zero, which results in the related parameters becoming close to zero. In figure \ref{simpleHierarchy}, for example, when the hyperparameters in the first level become small (compare to the hyperparameters in the second level), the model reduces to simple MNL. In contrast, when these hyperparameters are relatively large, the model reinforces our assumption of hierarchical classes. For this hierarchy we use the following priors:
\begin{eqnarray*}
\alpha_j | \tau_{0} \  & \sim & N(0, \tau_{0}^{2}) \ \ \qquad j = 1, ..., 4\\
{\phi}_{mkl} | \tau_{m} & \sim & N({0}, \tau_{m}^{2}) \\
\tau_{0}^{-2} & \sim & Gamma(a_{0}, b_{0})  \\
\tau_{m}^{-2} & \sim & Gamma(a_m, b_m)  \\
\end{eqnarray*} 
Here, $\phi_{mbl}$ refers to the parameter related to covariate $x_l$ and branch $k$ of node $m$. As we can see, all the branches that share the same node are controlled by one hyperparameter.

By introducing prior correlations between parameters for nearby classes, we can better handle situations in which these classes are hard to distinguish. If the hierarchy actually does provide information about how distinguishable classes are, we expect that performance will be improved. This would be especially true when the training set is small and the prior has relatively more influence on the results. Using an inappropriate hierarchy will likely lead to worse performance than a standard MNL model, but since the hyperparameters can adapt to reduce the prior correlations to near zero, the penalty may not be large. 

\section{\hspace*{-7pt}Results for synthetic datasets}\label{sec-many}\vspace{-12pt}

So far, we have discussed three alternative models: MNL, treeMNL, and corMNL. We first compare these models using a synthetic four-way classification problem with two covariates. Data are generated from each of these models in turn, and then fit with each model in order to test the robustness of the models when applied to data generated from other models . 

All regression parameters are given normal priors with mean zero. For the MNL model, the standard deviation for all the intercepts, $\tau_0$, and the standard deviation for the rest of coefficients, $\tau$, have the following priors:
\begin{eqnarray*}
\tau^{-2}_{0} & \sim & Gamma(1, 10)   \qquad (0.16, 0.38, 1.98) \\
\tau^{-2} & \sim & Gamma(1, 1) \; \; \qquad (0.52, 1.20, 6.27)
\end{eqnarray*}
We use the parameterization of the Gamma distribution in which Gamma$(a, b)$ has density $f(x|a, b)= [b^a \Gamma(a)]^{-1} x^{a-1} e^{-x/b}$, for which the mean is $ab$ and the standard deviation is $a^{1/2}b$. We gave the 2.5, 50 and 97.5 percentiles of $\tau$ in parenthesis.

For treeMNL and corMNL models, we assume that classes are arranged in a hierarchical form as shown in Figure \ref{simpleHierarchy}. This hierarchy implies that while it might be easy to distinguish between groups \{1, 2\} and \{3, 4\}, further separation of classes might not be as easy. As mentioned above, the treeMNL model for this hierarchy is comprised of three nested logistic models. These models are: $P(y \in \{1, 2\} | \alpha_1, \boldsymbol{\beta_1}, x)$, $P(y = 1 | \alpha_2, \boldsymbol{\beta_2}, x, y \in \{1, 2\})$ and $P(y = 3 | \alpha_3, \boldsymbol{\beta_3}, x, y \in \{3, 4\})$. For corMNL model, we use the priors discussed in section 2. The variance of regression parameters $\beta$'s in treeMNL and $\phi$'s in corMNL are regarded as hyperparameters. For these two models, one hyperparameter controls all the parameter emerging from the same node. These hyperparameters are given the following prior distributions: 
\begin{eqnarray*}
\tau_{0}^{-2} & \sim & Gamma(1, 10) \; \; \qquad (0.16, 0.38, 1.98) \\
\tau_{1}^{-2} & \sim & Gamma(1, 5)  \; \; \; \; \qquad (0.23, 0.54, 2.82) \\
\tau_{2}^{-2} & \sim & Gamma(1, 20)  \; \; \qquad (0.05, 0.12, 0.63) \\
\tau_{3}^{-2} & \sim & Gamma(1, 20)   \; \; \qquad (0.05, 0.12, 0.63) \\
\end{eqnarray*} 
Here, $\tau_1$, $\tau_2$ and $\tau_3$ correspond to nodes 1, 2, and 3 respectively (Figure \ref{simpleHierarchy}). These parameters have a narrower prior compared to $\tau$ in the MNL model. This is to account for the fact that the role of $\beta$ in the MNL model is played by more than one parameter in treeMNL and corMNL. Moreover, the regression parameters in the second level of hierarchy have a relatively smaller standard deviation $\tau$. As a result, these parameters tend to be smaller, making separation of class 1 from 2 and class 3 from 4 more difficult. 
 
We do three tests, in which we assume that each of the MNL, treeMNL and corMNL is the correct model. This allow us to see how robust each model is when data actually come from a somewhat different model. For each test, we sample a set of parameters from the prior distribution of the corresponding model. Pairs of data items ($x^{(i)}, y^{(i)}$) are generated by first drawing 10000 independent samples ${x}_1^{(i)}, {x}_2^{(i)}$ from the uniform($-5, 5$) distribution and then assigning each data item to one of the four possible classes. The assignment is either based on the multinomial model (for data generated from the MNL and corMNL) or based on successive logistic models (for data generated from the treeMNL). 

All three models are trained on the first 100 data items and tested on the remaining 9900 items. The regression coefficients were sampled from their posterior distribution using MCMC methods with single-variable slice sampling \cite{neal03}. At each iteration, we used the ``stepping out'' procedure to find the interval around the current point and used the ``shrinkage'' procedure for sampling from the interval. Since the hyperparameters were given conjugate priors, direct Gibbs sampling could be used for them. For all tests we drew 1000 samples from the posterior distributions. We discarded the initial 250 samples and used the rest for prediction. Performance is measured in terms of average log-probability and error rate on the test set. 

The above procedure was repeated 100 times. Each time, new regression parameters were sampled from priors and new pairs of data items were created based on the assumed models. The average results (over 100 iterations) are presented in Table \ref{simResults}. In this table, each column corresponds to the model used for generating the data and each row corresponds to the model used for building the classifier. As we can see, the diagonal elements have the best performance in each column. That is, the model whose functional form matches the data generation mechanism performs significantly better than the other two models (all \emph{p}-values based on average log-probability are less than 0.01 using a paired \emph{t}-test with $n=100$). Moreover, the results show that when the samples are generated according to the MNL model (i.e., classes are unrelated), corMNL has a significantly (\emph{p}-value $<$ 0.001) better performance compared to treeMNL. Reducing the size of training sets, $N$, to $50$ increases the gap between the average performance of models but provides the same conclusions. The conclusions also remain the same when we use different priors and ranges of covariates.
 
\begin{table} [t]
\begin{center}
\begin{tabular}{l | c |c|c|c|c|c|c|}
& \multicolumn{2}{|c|}{\bfseries Data from MNL} & \multicolumn{2}{|c|}{ \bfseries Data from treeMNL} &  \multicolumn{2}{|c|}{\bfseries Data from corMNL} \\
N=100 & AvgLogProb  & Error \% & AvgLogProb  &  Error \%\ & AvgLogProb  &  Error \%\ \\
\hline\hline
MNL method & {\bfseries-0.7958}  & \  {\bfseries 32.9} &  -0.8918 & \ 41.6 &  -0.9168  & \ 40.6 \\
\hline
treeMNL method & -0.8489  & \ 35.0 & {\bfseries -0.8770}  & {\bfseries 41.3} &  -0.9113  & \ 40.6   \\
\hline 
corMNL method & -0.7996  & \ 32.9 & -0.8797  & \ 41.4 &  {\bfseries -0.9075}  & {\bfseries 40.5}  \\
\hline
\end{tabular}
\end{center}

\label{simResults2}
\end{table}%

\begin{table}
\begin{center}
\begin{tabular}{l | c |c|c|c|c|c|c|}
& \multicolumn{2}{|c|}{\bfseries Data from MNL} & \multicolumn{2}{|c|}{ \bfseries Data from treeMNL} &  \multicolumn{2}{|c|}{\bfseries Data from corMNL} \\
 N= 50 & AvgLogProb  & Error \% & AvgLogProb  &  Error \%\ & AvgLogProb  &  Error \%\ \\
\hline\hline
MNL method & {\bfseries -0.8041}  & {\bfseries  33.6} &  -0.8915  & \ 42.4 &  -0.9248  & \ 41.4 \\
\hline
treeMNL method & -0.8506  & \ 35.1 & {\bfseries-0.8779}  & {\bfseries 42.0} &  -0.9227  & \ 41.5   \\
\hline 
corMNL method & -0.8086  & \ 33.7 & -0.8794  & \ 42.1 &  {\bfseries -0.9167}  & {\bfseries 41.3}  \\
\hline

\end{tabular}
\end{center}
\caption{Comparison of models on simulated data based on average log-probability using training sets of size $N = 100$ and $N = 50$. }
\label{simResults}

\vspace{14pt}

\begin{center}
\begin{tabular}{l | c |c|c|c|c|c|c|}
& \multicolumn{2}{|c|}{\bfseries Data from MNL} & \multicolumn{2}{|c|}{ \bfseries Data from treeMNL} &  \multicolumn{2}{|c|}{\bfseries Data from corMNL} \\
 N= 100 & AvgLogProb  & Error \% & AvgLogProb  &  Error \%\ & AvgLogProb  &  Error \%\ \\
\hline\hline
MNL method & {\bfseries -0.2539}  & {\bfseries  10.1} &  -0.3473  & \ 12.3 &  -0.3106  & \ 11.3 \\
\hline
treeMNL method & -0.6837  & \ 23.1 & -0.2898  & {\bfseries 9.7} &  -0.3614  & \ 12.1   \\
\hline 
corMNL method & -0.2910  & \ 10.3 & {\bfseries-0.2854}  & \ 9.9 &  {\bfseries -0.2841}  & {\bfseries 9.7}  \\
\hline

\end{tabular}
\end{center}
\caption{Comparison of models using a more complex hierarchy.  }
\label{simResults3}
\end{table}%

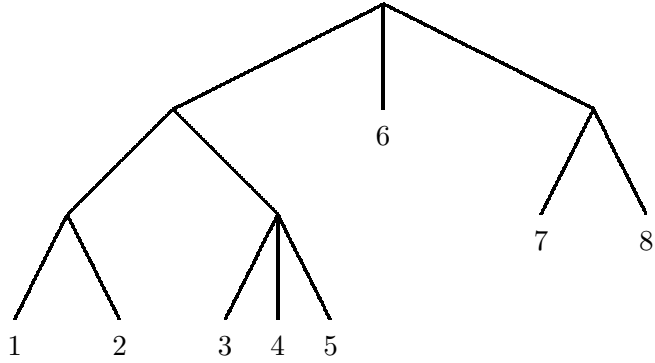
\begin{figure}[t]

\unitlength .7mm
\begin{picture}(120,70)(-45,0)
\linethickness{0.3mm}
\multiput(40,50)(0.24,0.12){167}{\line(1,0){0.24}}
\linethickness{0.3mm}
\multiput(20,30)(0.12,0.12){167}{\line(1,0){0.12}}
\linethickness{0.3mm}
\multiput(10,10)(0.12,0.24){83}{\line(0,1){0.24}}
\linethickness{0.3mm}
\multiput(20,30)(0.12,-0.24){83}{\line(0,-1){0.24}}
\linethickness{0.3mm}
\multiput(40,50)(0.12,-0.12){167}{\line(1,0){0.12}}
\linethickness{0.3mm}
\multiput(50,10)(0.12,0.24){83}{\line(0,1){0.24}}
\linethickness{0.3mm}
\put(60,10){\line(0,1){20}}
\linethickness{0.3mm}
\multiput(60,30)(0.12,-0.24){83}{\line(0,-1){0.24}}
\linethickness{0.3mm}
\put(80,50){\line(0,1){20}}
\linethickness{0.3mm}
\multiput(110,30)(0.12,0.24){83}{\line(0,1){0.24}}
\linethickness{0.3mm}
\multiput(120,50)(0.12,-0.24){83}{\line(0,-1){0.24}}
\linethickness{0.3mm}
\multiput(80,70)(0.24,-0.12){167}{\line(1,0){0.24}}
\put(10,5){\makebox(0,0)[cc]{1}}

\put(30,5){\makebox(0,0)[cc]{2}}

\put(50,5){\makebox(0,0)[cc]{3}}

\put(60,5){\makebox(0,0)[cc]{4}}

\put(70,5){\makebox(0,0)[cc]{5}}

\put(80,45){\makebox(0,0)[cc]{6}}

\put(110,25){\makebox(0,0)[cc]{7}}

\put(130,25){\makebox(0,0)[cc]{8}}

\end{picture}

\caption{A hypothetical hierarchy with a relatively more complex structure.}
\label{compHierarchy}
\end{figure}

While statistically significant, the results presented in Table \ref{simResults} might not be significant for some practical purposes. This is mostly due to the simplicity of the hierarchical structure. Next, we repeated the above tests with a more complex hierarchy, as shown in Figure \ref{compHierarchy}. For this problem we used four covariates randomly generated from the uniform(0, 1) distribution. In all three models, we used the same prior as before for the intercepts. For the MNL model we set $\tau^{-2} \sim Gamma(1, 1)$. The hyperparameters of treeMNL and corMNL were given $Gamma(1, 5)$, $Gamma(1, 20)$ and $Gamma(1, 100)$ priors for the first, second and third level of the hierarchy respectively.   

Table \ref{simResults3} shows the average results over 100 datasets for each test. As we can see, the differences between models are more accentuated. Nevertheless, corMNL is still performing well even when the data are generated by other models, where corMNL is outperformed only by the true model. When data come from treeMNL, the results from corMNL are very close to those of the actual model (i.e., treeMNL). 

\begin{figure}[b]
\unitlength .95mm
\begin{picture}(165,90)(-2,0)
\linethickness{0.15mm}
\multiput(65,70)(0.15,0.12){167}{\line(1,0){0.15}}
\linethickness{0.15mm}
\put(35,30){\line(0,1){20}}
\linethickness{0.15mm}
\put(155,60){\line(0,1){10}}
\put(5,45){\makebox(0,0)[cc]{Text}}

\put(33,65){\makebox(0,0)[cc]{Abstract}}

\put(51,25){\makebox(0,0)[cc]{Bullet}}

\put(51,20){\makebox(0,0)[cc]{Item}}

\put(40,55){\makebox(0,0)[cc]{Auth.}}

\put(40,50){\makebox(0,0)[cc]{List}}

\put(71,55){\makebox(0,0)[cc]{Subsec.}}

\put(71,50){\makebox(0,0)[cc]{Head.}}

\put(30,15){\makebox(0,0)[cc]{Fig.}}

\put(30,10){\makebox(0,0)[cc]{Cap.}}

\put(40,15){\makebox(0,0)[cc]{Table}}

\put(40,10){\makebox(0,0)[cc]{Cap.}}

\put(20,45){\makebox(0,0)[cc]{Ref.}}

\put(125,65){\makebox(0,0)[cc]{Main}}

\put(125,60){\makebox(0,0)[cc]{Title}}

\put(20,25){\makebox(0,0)[cc]{Foot}}

\put(20,20){\makebox(0,0)[cc]{Note}}

\put(60,55){\makebox(0,0)[cc]{Sec.}}

\put(60,50){\makebox(0,0)[cc]{Head.}}

\put(105,65){\makebox(0,0)[cc]{Eq.}}

\put(95,65){\makebox(0,0)[cc]{Eq.}}

\put(140,65){\makebox(0,0)[cc]{Decoration}}

\put(150,45){\makebox(0,0)[cc]{Graph}}

\put(162,45){\makebox(0,0)[cc]{Fig.}}

\put(145,55){\makebox(0,0)[cc]{Table}}

\put(164,55){\makebox(0,0)[cc]{Code}}

\linethickness{0.3mm}
\multiput(40,60)(0.12,0.24){42}{\line(0,1){0.24}}
\linethickness{0.3mm}
\multiput(45,70)(0.12,-0.24){42}{\line(0,-1){0.24}}
\linethickness{0.3mm}
\multiput(20,70)(0.42,0.12){167}{\line(1,0){0.42}}
\linethickness{0.3mm}
\multiput(5,50)(0.12,0.16){125}{\line(0,1){0.16}}
\linethickness{0.3mm}
\multiput(20,70)(0.12,-0.16){125}{\line(0,-1){0.16}}
\linethickness{0.3mm}
\multiput(155,70)(0.12,-0.12){83}{\line(1,0){0.12}}
\linethickness{0.3mm}
\multiput(150,50)(0.12,0.24){42}{\line(0,1){0.24}}
\put(50,55){\makebox(0,0)[cc]{Ed.}}

\put(50,50){\makebox(0,0)[cc]{List}}

\put(30,35){\makebox(0,0)[cc]{}}

\put(15,50){\makebox(0,0)[cc]{}}

\put(115,65){\makebox(0,0)[cc]{Page}}

\put(90,80){\makebox(0,0)[cc]{}}

\linethickness{0.3mm}
\multiput(75,70)(0.12,0.16){125}{\line(0,1){0.16}}
\linethickness{0.3mm}
\multiput(20,30)(0.12,0.16){125}{\line(0,1){0.16}}
\linethickness{0.3mm}
\multiput(35,50)(0.12,-0.16){125}{\line(0,-1){0.16}}
\put(55,65){\makebox(0,0)[cc]{Header}}

\put(75,65){\makebox(0,0)[cc]{Footer}}

\put(35,30){\makebox(0,0)[cc]{}}

\linethickness{0.3mm}
\multiput(60,60)(0.12,0.24){42}{\line(0,1){0.24}}
\linethickness{0.3mm}
\multiput(65,70)(0.12,-0.24){42}{\line(0,-1){0.24}}
\put(85,50){\makebox(0,0)[cc]{}}

\put(82,55){\makebox(0,0)[cc]{Fig.}}

\put(93,55){\makebox(0,0)[cc]{Table}}

\put(82,50){\makebox(0,0)[cc]{Label}}

\put(93,50){\makebox(0,0)[cc]{Label}}

\linethickness{0.3mm}
\multiput(30,20)(0.12,0.24){42}{\line(0,1){0.24}}
\linethickness{0.3mm}
\multiput(35,30)(0.12,-0.24){42}{\line(0,-1){0.24}}
\linethickness{0.3mm}
\multiput(145,60)(0.12,0.12){83}{\line(1,0){0.12}}
\linethickness{0.3mm}
\multiput(155,60)(0.12,-0.24){42}{\line(0,-1){0.24}}
\linethickness{0.3mm}
\put(20,50){\line(0,1){20}}
\linethickness{0.3mm}
\multiput(85,70)(0.12,0.48){42}{\line(0,1){0.48}}
\linethickness{0.3mm}
\multiput(35,70)(0.33,0.12){167}{\line(1,0){0.33}}
\linethickness{0.3mm}
\multiput(90,90)(0.15,-0.12){167}{\line(1,0){0.15}}
\linethickness{0.3mm}
\multiput(90,90)(0.21,-0.12){167}{\line(1,0){0.21}}
\linethickness{0.3mm}
\multiput(45,70)(0.27,0.12){167}{\line(1,0){0.27}}
\linethickness{0.3mm}
\multiput(90,90)(0.39,-0.12){167}{\line(1,0){0.39}}
\linethickness{0.3mm}
\multiput(90,90)(0.3,-0.12){167}{\line(1,0){0.3}}
\linethickness{0.3mm}
\multiput(80,60)(0.12,0.24){42}{\line(0,1){0.24}}
\linethickness{0.3mm}
\multiput(85,70)(0.12,-0.24){42}{\line(0,-1){0.24}}
\linethickness{0.3mm}
\multiput(55,70)(0.21,0.12){167}{\line(1,0){0.21}}
\linethickness{0.3mm}
\multiput(90,90)(0.12,-0.48){42}{\line(0,-1){0.48}}
\linethickness{0.3mm}
\multiput(90,90)(0.12,-0.16){125}{\line(0,-1){0.16}}
\put(105,60){\makebox(0,0)[cc]{\#}}

\put(115,60){\makebox(0,0)[cc]{\#}}

\end{picture}

\caption{Hierarchical structure of document labels.}
\label{documentLabelling}
\end{figure}
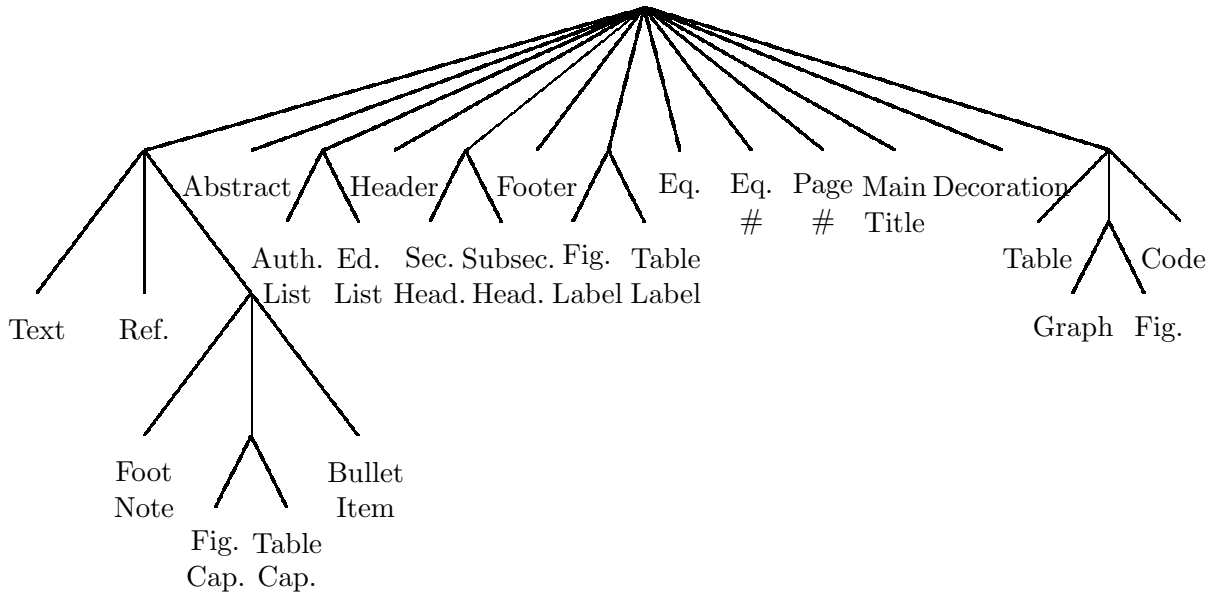

\section{\hspace*{-7pt}Results on a document labelling problem}\label{sec-tests}\vspace{-12pt}

We next test our approach using a dataset of the page images of 15 articles (472 pages) from the Journal of Machine Learning Research (JMLR) collected by Laven \cite{laven05}. Each page was segmented to several regions, and each region was manually classified to one of 24 possible classes. Figure  \ref{documentLabelling} presents these classes in a hierarchical form. The hierarchy is based on our belief regarding how difficult it is to separate classes from each other using the available features.   

The covariates are 59 different features such as the location of the region on the page and the density of the ink inside the region. We normalized all features so they have zero mean and standard deviation 1. 

Laven divided the dataset into a training set (including 10 articles with 3521 regions), and a test set (including 5 articles with 2035 regions). The items from the same article are considered as independent even though they clearly are not. Although this may cause overfitting problems, we follow the same assumption for our initial test in order to make our results comparable to Laven's. 

We trained the three models (MNL, corMNL and treeMNL) on the training set and evaluated their performance on the test set. The coefficient parameters in the MNL models were given normal prior with mean zero. The variances of these parameters were regarded as hyperparameters. For this problem, since the number of covariates, $p=59$, is relatively large, we use the Automatic Relevance Determination (ARD) method suggested by \cite{neal96}. ARD employs a hierarchical prior to determine how relevant each covariate is in classification of objects. In this method, one hyperparameter, $\sigma_l$, is used to control the variance of all coefficients, $\beta_{jl}$ ($j=1, ..., c$), for covariate $x_l$. If a covariate is irrelevant, its hyperparameter will tend to be small, forcing the coefficients for that covariate be near zero. As before, we also use one hyperparameter, $\tau$, to control all $\beta$'s in the MNL model. We set the standard deviation of $\beta_{jl}$ equal to $\tau \sigma_{l}$. Therefore, while $\sigma_l$ control the relevance of covariate $x_l$ compared to other covariates, the scale parameter $\tau$, controls the overall usefulness of all covariates in separating classes. For the MNL model we used the following priors:
\begin{eqnarray*}
\alpha_j | \; \ \; \sigma_{0} & \sim & N(0, \tau_{0}^{2}) \quad \qquad \qquad j = 1, ..., 24\\
\beta_{jl} | \tau, \sigma_l & \sim & N({0}, \tau^{2}\sigma_{l}^{2})  \qquad  \qquad l  = 1, ..., 59\\
\tau^{-2}_{0} & \sim & Gamma(0.5, 1)  \qquad (0.63, 2.09, 46.31)\\
\tau^{-2} & \sim & Gamma(0.5, 20)  \; \ \quad (0.14, 0.47, 10.07)\\
\sigma^{-2}_{l} & \sim & Gamma(1, 10) \  \qquad (0.16, 0.38, 1.98)
\end{eqnarray*} 

Similar priors are used for the parameters of treeMNL and corMNL. For these two models, we used one hyperparamter, $\sigma^{-2}_{l} \sim Gamma(1, 10)$ to control all parameters related to covariate $x_l$. We also used one scale parameter $\tau^{-2}_{m} \sim Gamma(0.5, 100)$ for all parameters ($\beta$'s in treeMNL, $\phi$'s in corMNL) sharing the same node $m$. The prior for the intercepts was the same as in the MNL model. 

We used Hamiltonian dynamics \cite{neal93} for sampling from the posterior distribution of coefficients parameters. To reduce the random walk aspect of sampling procedure, we use a reasonably large number of leapfrog steps ($L=500$). The stepsizes is set to 0.02 in order to maintain an acceptance rate of about 90\%. In the MNL and corMNL models, the new updates are proposed for all regression parameters simultaneously. Nested MNL models in treeMNL are updated separately since they are regarded as independent models. The coefficient parameters within each nested model, however, are updated at the same time. We used single-variable slice sampling \cite{neal03}, as described above, for hyperparameters. The convergence of the Markov chain simulations were assessed by plotting the values of hyperparameters and average log-likelihood (on training cases). We ran each chain for 2500 iterations, of which the first 500 were discarded.

\begin{table} [h]
\begin{center}
\begin{tabular}{l | c |c|c|}
&  AvgLogProb &  Error rate \% \\
\hline\hline
Laven's model (ML) & $-$ & 9.4\\
\hline
MNL &  -0.3998 & 9.4 \\
\hline
treeMNL & -0.3602 & 9.3 \\
\hline
corMNL & \textbf{-0.3598} & \textbf{9.2} \\
\hline
\end{tabular}
\end{center}
\caption{Performance of models for the document labelling problem on the test set.}
\label{finalRes1}
\end{table}%

Table \ref{finalRes1} compares the results from different models. As we can see, our MNL model has exactly the same performance as the likelihood-based model developed by Laven. The corMNL and treeMNL models have a slightly better error rate and a higher average log-probability compared to the MNL model.    

We speculated that the assumed structure would have a more profound effect when a smaller sample is available. To test this idea, we reduced the size of the training set to 200, and used 10 non-overlapping training sets in order to obtain a confidence interval for the performance of each model. To avoid problems from incorrectly assuming independence, we first combined the training and the test sets and then randomly sampled (without replacement) 10 training sets each with 200 cases and used the remaining 3556 cases as the test set. Sampling was performed without regard to which document each region comes from, which artificially makes the independence assumption be true. We trained the models on each training set and measured the performance on the test set. We discarded a burn-in period of 500 initial updates and retained the next 3500 updates for prediction. The results are shown in Table \ref{finalRes2}. Since the test set is large, the observed variation in the results is dominated by the difference between training sets, and we can use paired $t$-tests ($n=10$) for comparing the results. Based on the average log-probability results, the corMNL model outperforms both MNL (\emph{p}-value $=$ 0.002 ) and treeMNL (\emph{p}-value $=$ 0.004) models. The improvement on error rate is only slightly significant compared to MNL (\emph{p}-value = 0.09) and treeMNL (\emph{p}-value = 0.06). The differences in performance of MNL and treeMNL are not statistically significant.
\begin{table} [h]
\begin{center}
\begin{tabular}{l | c |c|c|}
& AvgLogProb & Error rate \% \\
\hline\hline
MNL &  -0.5441  &  15.4 \\
\hline
treeMNL & -0.5381  & 15.5  \\
\hline
corMNL & \textbf{-0.5160 } & \textbf{14.9 } \\
\hline
\end{tabular}
\end{center}
\caption{Performance of models for the document labelling problem on the test set using 10 non-overlapping training sets of size 200.}
\label{finalRes2}
\end{table}%

\section{\hspace*{-7pt}Conclusions and future directions}\label{sec-conc}\vspace{-12pt}
In this paper, we have introduced a new approach for modelling hierarchical classes. Our analysis shows that when the hierarchy actually does provide information regarding the similarity of classes, our approach outperforms simple MNL model and models based on decomposing the hierarchy into nested MNL models. Our method can be applied to many classification problems where there is prior knowledge regarding the structure of classes. One such problem, which we hope to investigate soon, is annotation of gene functions. 

So far, we have focused only on simple tree-like structures. There are other hierarchical structures that are more complex than a tree. For example, one of the most commonly used gene annotation schemes, known as Gene Ontology (GO), is implemented as a directed acyclic graph (DAG). In this structure a node can have more than one parent. Our method, as it is, cannot be applied to these problems, but it should be possible to extend the idea of summing coefficients along a path to the class in order to allow for multiple paths. 
	
In our approach, we considered only one structure for each hierarchical classification problem. However, we might sometimes be able to think of more than one possible class hierarchy. We might then consider all possible hierarchies as equally likely and form an ensemble model based on this assumption (as done in \cite{eibe04}), but this might not be easy when the number of classes is large. Also, in most cases we have prior knowledge that leads us to prefer some structures over others a priori. It is possible to generalize our method to multiple hierarchies. As for the generalization to DAG's, it should be possible to sum coefficients along the multiple paths within different hierarchies. We can further use a set of hyperparameters to discover the relevance of each hierarchy.  

The results presented in this paper are for linear models. We expect that a similar approach can be used in non-linear models such as neural networks. 

\section*{Acknowledgements}\vspace{-12pt}
We thank Kevin Laven for providing the document labelling dataset. This research was supported by the Natural Sciences and Engineering Research Council of Canada. Radford Neal holds a Canada Research Chair in Statistics and Machine Learning.

\end{document}